\documentclass[11pt]{amsart}

\usepackage{amsmath,amssymb,amsthm}
\usepackage{mathtools}
\usepackage{booktabs}
\usepackage{graphicx}
\usepackage{subcaption}
\usepackage[margin=1in]{geometry}

\newtheorem{theorem}{Theorem}[section]
\newtheorem{lemma}[theorem]{Lemma}
\newtheorem{proposition}[theorem]{Proposition}

\newcommand{\RR}{\mathbb{R}}
\newcommand{\E}{\mathbb{E}}
\newcommand{\sig}{\sigma}
\newcommand{\norm}[1]{\left\|#1\right\|}

\newcommand{\xt}{\mathbf{x}_{*}}
\newcommand{\x}{\mathbf{x}}
\newcommand{\bb}{\mathbf{b}}
\newcommand{\aaa}{\mathbf{a}}

\begin{document}

\title{Two-point residual step size gradient methods for linear systems}

\author{Xin Xu}
\address{School of Mathematics and Statistics, Shanxi University, No.~92 Wucheng Road, Xiaodian District, Taiyuan 030006, Shanxi, China}
\email{xustonexin@163.com}
\date{\today}

\begin{abstract}
Kaczmarz-type methods often converge slowly when the coefficient matrix is ill-conditioned. In particular, when the smallest singular value is small, every row of the coefficient matrix is nearly orthogonal to the corresponding right singular direction, so that the iteration stalls along that direction. To break this slow-convergence regime, we borrow the delayed-step-size idea of the Barzilai-Borwein (BB) method and, within the framework of the residual-based surrogate-hyperplane Kaczmarz method, propose two new step sizes for gradient-descent-type methods solving linear systems: a delayed residual step size \(\alpha_k^{\mathrm{DR}}\) and a residual-difference step size \(\alpha_k^{\mathrm{RD}}\), which correspond, respectively, to the BB1 and BB2 step sizes in convex quadratic optimization. We prove that \(\alpha_k^{\mathrm{RD}}\le\alpha_k^{\mathrm{DR}}\), and we analyze the convergence of the resulting method in a canonical model with two singular values: although the step sizes are non-monotone, the residual norm is shown to decrease after finitely many steps. Numerical experiments on some ill-conditioned linear systems illustrate the effectiveness of the proposed methods.
\end{abstract}


\maketitle

\section{Introduction}\label{sec:intro}

We consider the solution of the consistent linear system 
\begin{equation}\label{eq:linearsys}
A\x=\bb, 
\end{equation}
where $A\in\RR^{m\times n}$ has full rank, $\x\in\RR^{n}$, and $\bb\in\RR^{m}$. Let $\xt=A^{\dagger}\bb$ be the unique least-norm solution, where $A^{\dagger}$ denotes the Moore-Penrose pseudoinverse of $A$. Numerous problems in science and engineering can be formulated as the solution of a linear system, such as the discretization of partial differential equations \cite{saad2011}, image reconstruction \cite{gordon1970,han2021}, phase retrieval \cite{chi2016}, and large-scale least-squares data fitting \cite{needell2014}. Since such systems are frequently too large for direct solvers, the development of efficient iterative solvers is of considerable practical importance. As an effective iterative scheme, the randomized Kaczmarz (RK)~\cite{strohmer2009} algorithm is widely used for this purpose. Let $\aaa_{i_k}^{\top}$ denote the $i_k$-th row of $A$. At the $k$-th step, the RK algorithm selects a row index $i_k$ and projects the current approximation $\x_k$ orthogonally onto the hyperplane defined by that row, yielding the update
\begin{equation}
\x_{k+1}=\x_k+\omega_{k}\frac{\bb_{i_k}-\aaa_{i_k}^{\top}\x_k}{\norm{\aaa_{i_k}}^2}\aaa_{i_k},
\end{equation}
where $\omega_{k}\in(0,2)$. When the relaxation parameter satisfies $\omega_{k}=1$ and the rows are visited cyclically, i.e. $i_k=(k~\text{mod}~m)+1$, the scheme reduces to the original Kaczmarz (KACZ)~\cite{kacz1937} algorithm; when $i_k$ is drawn with probability proportional to $\norm{\aaa_{i_k}}^2$, one obtains the RK algorithm, and its expected contraction is
\begin{equation}\label{eq:rkr}
\E\,\norm{\x_{k+1}-\xt}^2\;\le\;\Big(1-\frac{\sig_{\min}(A)^2}{\norm{A}_F^2}\Big)\norm{\x_k-\xt}^2,
\end{equation}
so the convergence is governed by the ratio $\sig_n^2/\norm{A}_F^2$. When $A$ is ill-conditioned, i.e.\
$\sig_n\ll\sig_1$, the method is extremely slow, where $\sig_n$ and $\sig_1$ are the smallest and largest singular values of $A$, respectively.

Owing to the simplicity, the transparent geometric interpretation, and the favorable numerical performance of such projection-type iterative methods, this family of algorithms has attracted considerable attention in recent years, and numerous block, greedy, and momentum-accelerated variants have been proposed to improve their efficiency and to broaden their range of applicability; see, e.g., \cite{gower2015,needell2014,du2020,bai2018,liu2016} and the references therein.

In order to exploit the information contained in $A$ as fully as possible, Wang and Yin \cite{Wang2023} proposed a residual-based surrogate-hyperplane Kaczmarz (RSHK) method, whose iteration takes the form 
\begin{equation}\label{equ:RSHK}
	\x_{k+1}=\x_{k}+\frac{\mathbf{r}_{k}^{\top}(\bb-A\x_{k})}{\|A^{\top}\mathbf{r}_{k}\|_{2}^{2}}A^{\top}\mathbf{r}_{k}=\x_{k}+\frac{\mathbf{r}_{k}^{\top}\mathbf{r}_{k}}{\|A^{\top}\mathbf{r}_{k}\|_{2}^{2}}A^{\top}\mathbf{r}_{k},
\end{equation}
where the coefficient
\begin{equation}\label{equ:rrk}
\mathbf{r}_{k}=\bb-A\x_{k}\in\RR^{m}
\end{equation}
is the residual vector; that is, each step requires the solution of the compressed linear equation
\begin{equation}\label{equ:reducedLinearSys}
	\mathbf{r}_{k}^{\top}Ax=\mathbf{r}_{k}^{\top}\bb.
\end{equation}
Note that replacing the coefficient vector $\mathbf{r}_{k}$ with different vectors gives rise to different surrogate-hyperplane KACZ algorithms. For instance, when $\mathbf{r}_{k}$ is chosen as the unit vector $\mathbf{e}_{i_{k}}$ (i.e., the $i_k$-th column of the identity matrix $I$), the algorithm reduces to the (randomized) KACZ method; further choices can be found in \cite{Wang2023} and the references therein. A closer examination reveals that the method weights all rows of $A$ by the residual, thereby producing a surrogate hyperplane, and it therefore belongs to the class of residual-adaptive iterative algorithms.

As pointed out in \cite{Wang2023}, an important observation is that, from the viewpoint of gradient descent, the gradient direction of the objective function $\frac{1}{2}\|A\x-\bb\|_{2}^{2}$ at iterate $\x_k$ is given by
\begin{equation}\label{equ:gk}
	\mathbf{g}_{k}=A^{\top}(A\x-\bb)=-A^{\top}\mathbf{r}_{k}.
\end{equation}
Setting the step size to $\alpha_{k}$, the gradient descent iteration takes the form 
\begin{equation}\label{equ:iterate}
	\x_{k+1}=\x_{k}-\alpha_{k}\mathbf{g}_{k}=\x_{k}+\alpha_{k}A^{\top}\mathbf{r}_{k}.
\end{equation}
To determine the step size $\alpha_{k}$, we minimize
\begin{equation*}
	\|\x_{k+1}-\x_{*}\|_{2}^{2}=\|\x_{k}-\x_{*}+\alpha_{k}A^{\top}\mathbf{r}_{k}\|_{2}^{2}=\|\x_{k}-\x_{*}\|_{2}^{2}+\alpha_{k}^2\mathbf{r}_{k}^{\top}AA^{\top}\mathbf{r}_{k}-2\alpha_{k}\mathbf{r}_{k}^{\top}\mathbf{r}_{k}
\end{equation*}
with respect to $\alpha_{k}$, which yields 
\begin{equation}\label{equ:SD}
	\alpha_{k}^{\mathrm{RSHK}}=\frac{\mathbf{r}_{k}^{\top}\mathbf{r}_{k}}{\|A^{\top}\mathbf{r}_{k}\|_{2}^{2}}
\end{equation}
from the first-order optimality condition.
Therefore, the RSHK algorithm corresponds precisely to the method of steepest descent (SD) within the class of gradient descent methods. This observation offers a valuable clue for solving linear systems from an optimization perspective.

Achiya Dax~\cite{Achiya2023} revealed a surprising, even counterintuitive, feature of the Kaczmarz method, namely the so-called Kaczmarz's anomaly: when the number of equations $m$ is much smaller than the number of unknowns $n$, the method converges rapidly; as $m$ increases towards $n$, convergence slows down and, in the limit $m\to n$, the number of iterations grows sharply; yet once $m$ exceeds $n$, entering the overdetermined regime, convergence accelerates with the number of equations, and rapid convergence is recovered when $m\gg n$. Since this behavior runs counter to intuition, it has been termed Kaczmarz's anomaly.

To explain this phenomenon, \cite{Achiya2023} proves, by the interlacing property of singular values, that as the number of rows $i$ approaches $n$ from below, the smallest singular value of $A_i$ decreases sharply (and increases again once $i$ passes $n$), while the associated condition number rises sharply and then falls rapidly—the so-called smallest singular value anomaly and condition number anomaly. Subsequently, by connecting the Kaczmarz method with the associated SOR (successive overrelaxation) method, \cite{Achiya2023} derives an upper bound on the spectral radius of the iteration matrix, showing that a small condition number yields fast convergence, and further shows that a small value of the smallest singular value entails slow convergence. In this way, the anomaly is accounted for.

This naturally raises a fundamental question: why does the smallest singular value exert such a significant influence on the convergence of KACZ-type algorithms?

The present paper attempts to answer this question, and further shows that the anomaly of the KACZ method stems from an inherent feature of the algorithm itself. Before describing our approach, we recall the Barzilai-Borwein (BB) method \cite{barzilai1988}, a well-known acceleration of steepest descent for unconstrained minimization. Whereas the method of steepest descent determines the step size by an exact line search, the BB method determines it by enforcing a quasi-Newton (secant) condition between the two most recent iterates; for the quadratic $f(\mathbf{x})=\tfrac{1}{2}\mathbf{x}^{\top}B\mathbf{x}-\mathbf{b}^{\top}\mathbf{x}$, the two BB step sizes are given by the reciprocals of the corresponding Rayleigh quotients
\begin{equation}\label{equ:BB}
\alpha_k^{\mathrm{BB1}}=\frac{\mathbf{s}_{k-1}^{\top}\mathbf{s}_{k-1}}{\mathbf{s}_{k-1}^{\top}\mathbf{y}_{k-1}}=\frac{\mathbf{g}_{k-1}^{\top}\mathbf{g}_{k-1}}{\mathbf{g}_{k-1}^{\top}B\mathbf{g}_{k-1}},\qquad
\alpha_k^{\mathrm{BB2}}=\frac{\mathbf{s}_{k-1}^{\top}\mathbf{y}_{k-1}}{\mathbf{y}_{k-1}^{\top}\mathbf{y}_{k-1}}=\frac{\mathbf{g}_{k-1}^{\top}B\mathbf{g}_{k-1}}{\mathbf{g}_{k-1}^{\top}B^{2}\mathbf{g}_{k-1}},	
\end{equation}
where 
\begin{equation}\label{equ:sy}
\mathbf{g}_{k}=\nabla f(\x_{k}),\quad  \mathbf{s}_{k-1}=\mathbf{x}_{k}-\mathbf{x}_{k-1},\quad \text{and}\quad \mathbf{y}_{k-1}=\mathbf{g}_{k}-\mathbf{g}_{k-1}.	
\end{equation}
A salient feature of the BB method is its non-monotonicity: for the sequence of iterates it generates, neither the gradient norms nor the objective function values are guaranteed to decrease monotonically, and transient increases are permitted. Such transient increases, however, are typically accompanied by substantial subsequent decreases in the gradient norms or the objective function values. Overall, these non-monotone step sizes can dramatically accelerate convergence relative to the SD step size
\begin{equation}\label{equ:SDD}
\alpha_{k}^{SD}=\frac{\mathbf{g}_{k}^{\top}\mathbf{g}_{k}}{\mathbf{g}_{k}^{\top}A\mathbf{g}_{k}},	
\end{equation}
especially on some ill-conditioned quadratics, and the method is known to converge $R$-linearly on strongly convex quadratics \cite{dai2002}. For research on the BB method and its extensions, see \cite{Raydan1997BarzilaiBorweinGradienta,Zhou2006GradientMethodsAdaptive,Huang2021EquippingBarzilaiBorwein,An2025RegularizedBarzilaiBorwein,Xu2025PBB} and the references therein. 

The success of the BB step size naturally raises the question of whether the idea underlying the BB method can be applied to accelerate the Kaczmarz algorithm. To this end, we draw upon the optimization viewpoint of the RSHK method and the spirit of the BB method to propose two new iterative algorithms, which aim to partially remedy the anomalous behavior of the Kaczmarz method.

The remainder of this paper is organized as follows. In Section 2, we analyze the slow convergence mechanism of the Kaczmarz method. In Section 3, we formally propose two new algorithms and provide a convergence analysis for a simple case. In Section 4, we conduct numerical comparative experiments to verify the effectiveness of the two proposed algorithms. In Section 5, we conclude the paper and give an outlook for future work.

\section{Slow convergence of the Kaczmarz method} 
Although Kaczmarz-type methods are effective for solving general linear systems, their performance is strongly influenced by the condition number of the coefficient matrix, and in particular by its smallest singular value. In this section, we analyze why Kaczmarz-type methods converge slowly when the smallest singular value is small.

\begin{proposition}\label{prop:orth}
	Let $A\in\mathbb{R}^{m\times n}$ with $m\ge n$ and $\mathrm{rank}(A)=n$, and let 
	\begin{equation}\label{equ:svd}
	A=U\Sigma V^T,
	\end{equation}
	be the economy singular value decomposition of $A$, where
	\begin{equation*}
	\Sigma=\mathrm{diag}(\sig_1,\dots,\sig_n)\in\RR^{n\times n},\quad U=[\mathbf{u}_1,\dots,\mathbf{u}_n]\in\RR^{m\times n},\quad
	V=[\mathbf{v}_1,\dots,\mathbf{v}_n]\in\RR^{n\times n},	
	\end{equation*}
	with $\sig_1\ge\cdots\ge\sig_n>0$, $U^{\top}U=I_n$, $V^{\top}V=VV^{\top}=I_n$. 

	$\mathrm{(i)}$ For all $1\le i\le m$ and $1\le j\le n$,
	$\mathbf{a}_i^T\mathbf{v}_j=\sigma_j\,\mathbf{u}_{j}^{(i)}$, $\big|\mathbf{a}_i^T\mathbf{v}_j\big|\le\sigma_j.$
	
	$\mathrm{(ii)}$ Let $1\le k\le n-1$ and write
	$S_k:=\mathrm{span}\{\mathbf{v}_{n-k+1},\dots,\mathbf{v}_n\}$ and
	$S_k^\perp:=\mathrm{span}\{\mathbf{v}_1,\dots,\mathbf{v}_{n-k}\}$. Then for every $i=1,\dots,m$,
	$\mathrm{dist}\big(\mathbf{a}_i,\,S_k^\perp\big)\le\sigma_{n-k+1}.$
	In particular, when the trailing $k$ singular values $\sigma_{n-k+1}\ge\cdots\ge\sigma_n$ are all small, every row of $A$ is nearly orthogonal to the subspace $S_k$.
\end{proposition}
\begin{proof}
	Since $A\mathbf{v}_j=\sigma_j\mathbf{u}_j$ and $\mathbf{a}_i^T=\mathbf{e}_i^T A$, we have $\mathbf{a}_i^T\mathbf{v}_j=\sigma_j\,\mathbf{u}_{j}^{(i)}$, which proves (i). Because $\{\mathbf{v}_j\}$ is an orthonormal basis of $\mathbb{R}^n$,
	$\mathbf{a}_i=\sum_{j=1}^{n}\sigma_j\,\mathbf{u}_{j}^{(i)}\,\mathbf{v}_j,$
	the orthogonal projection of $\mathbf{a}_i$ onto $S_k^\perp$ is
	$\sum_{j=1}^{n-k}\sigma_j\,\mathbf{u}_{j}^{(i)}\mathbf{v}_j$. Hence
	$\mathrm{dist}\big(\mathbf{a}_i,\,S_k^\perp\big)^2
	=\sum_{j=n-k+1}^{n}\sigma_j^2\,\big(\mathbf{u}_{j}^{(i)}\big)^2
	\le\sigma_{n-k+1}^2\sum_{j=n-k+1}^{n}\big(\mathbf{u}_{j}^{(i)}\big)^2
	\le\sigma_{n-k+1}^2,$
	where the last inequality uses $U^{\top}U=I_n$, which yields $\sum_{j=1}^{n}(\mathbf{u}_{j}^{(i)})^2\le1$.
\end{proof}

Observation. As the trailing $k$ singular values tend to zero, Proposition~\ref{prop:orth}(ii) shows that every row of $A$ becomes nearly orthogonal to $S_k=\mathrm{span}\{\mathbf{v}_{n-k+1},\dots,\mathbf{v}_n\}$. Since a Kaczmarz-type update moves only along the rows (or their linear combinations), the progress in the directions $\{\mathbf{v}_{n-k+1},\dots,\mathbf{v}_n\}$ is correspondingly small; this suggests that the iteration stalls along these directions, which explains its slow convergence.

The influence of the smallest singular value on Kaczmarz-type algorithms described above is, in essence, inherent to the structure of the algorithms themselves: as long as the iterates are obtained by orthogonally projecting onto the hyperplanes corresponding to the rows of $A$, or onto certain linear combinations of these rows, this slow-convergence regime cannot be readily broken.

\section{New iterative algorithms}
In this section, we employ gradient-descent algorithms to break the slow-convergence regime of Kaczmarz-type methods. 

\subsection{Derive new step size}
As established in the preceding analysis, the RSHK method is equivalent to the steepest descent method. From the viewpoint of optimization, $A^{T}\mathbf{r}_k$ is the locally optimal descent direction, and the step size $\alpha_k^{\mathrm{RSHK}}$ \eqref{equ:SD} is the optimal one. As is well known from steepest descent in convex quadratic programming, however, the strategy of combining an optimal search direction with an optimal step size is generally not optimal in a global sense. Recall that the search direction of the BB method and that of the steepest descent method are both the negative gradient direction; yet the former step size is exactly the latter delayed by one step, namely $\alpha_k^{\mathrm{BB1}}=\alpha_{k-1}^{\mathrm{SD}}$. This simple delayed-step strategy elevates the first-order difference dynamical system associated with the steepest descent method to second order, thereby substantially improving the convergence performance \cite{Xu2025OnRBB}. Inspired by the delayed-step strategy of the BB method, we consider, on the basis of RSHK, the following two step sizes:
\begin{equation}\label{equ:RDDR}
\alpha_k^{\mathrm{DR}}=\frac{\mathbf{r}_{k-1}^{T}\mathbf{r}_{k-1}}{\|A^{T}\mathbf{r}_{k-1}\|_2^{2}}=\alpha_{k-1}^{\mathrm{RSHK}},\quad
\alpha_k^{\mathrm{RD}}=\frac{\mathbf{d}_{k-1}^{T}\mathbf{d}_{k-1}}{\|A^{T}\mathbf{d}_{k-1}\|_2^{2}},
\end{equation}
where $\mathbf{d}_{k-1}=\mathbf{r}_{k}-\mathbf{r}_{k-1}$. Here, DR and RD denote, respectively, the delayed residual and the difference of residuals. Formally, we denote the gradient descent method using \(A^{T}\mathbf{r}_{k}\) as the search direction and step size \(\alpha_{k}^{DR}\) as the DR method, and the counterpart with step size \(\alpha_{k}^{RD}\) as the RD method. It is worth emphasizing that the step size $\alpha_k^{\mathrm{DR}}$ is obtained directly by delaying the RSHK step size by one step. To see that $\alpha_k^{\mathrm{RD}}$ is the
counterpart of the BB2 step size in the present setting, note that, by \eqref{equ:gk},
\eqref{equ:BB} and \eqref{equ:sy},
\begin{equation*}
\mathbf{d}_{k-1}=\mathbf{r}_k-\mathbf{r}_{k-1}=-A\mathbf{s}_{k-1},\quad
\mathbf{y}_{k-1}=A^{T}(\mathbf{r}_{k-1}-\mathbf{r}_k)=-A^{T}\mathbf{d}_{k-1}=A^{T}A\mathbf{s}_{k-1}.
\end{equation*}
Consequently,
\begin{equation*} \mathbf{s}_{k-1}^{T}\mathbf{y}_{k-1}=\|A\mathbf{s}_{k-1}\|_2^{2}=\|\mathbf{d}_{k-1}\|_2^{2} \quad\text{and}\quad
\mathbf{y}_{k-1}^{T}\mathbf{y}_{k-1}=\|A^{T}\mathbf{d}_{k-1}\|_2^{2}, 
\end{equation*}
and hence
\begin{equation*}
\alpha_k^{\mathrm{RD}}=\frac{\mathbf{s}_{k-1}^{T}\mathbf{y}_{k-1}}{\mathbf{y}_{k-1}^{T}\mathbf{y}_{k-1}}
=\frac{\mathbf{d}_{k-1}^{T}\mathbf{d}_{k-1}}{\|A^{T}\mathbf{d}_{k-1}\|_2^{2}}.
\end{equation*}
This is precisely the expression given in \eqref{equ:RDDR}.

\begin{theorem}
	Let $A\in\mathbb{R}^{m\times n}$ with $m\ge n$ and $\operatorname{rank}(A)=n$, and assume that
	$A^{T}\mathbf{r}_{k-1}\neq\mathbf{0}$, so that the step sizes in \eqref{equ:RDDR} are well defined.
	Then
	\begin{equation}\label{equ:inequality}
	\alpha_{k}^{\mathrm{RD}}\le\alpha_{k}^{\mathrm{DR}}.
	\end{equation}
\end{theorem}

\begin{proof}
	It follows from \eqref{equ:iterate} that
	\begin{equation}\label{equ:rk}
	\mathbf{r}_{k+1}=\mathbf{r}_{k}-\alpha_{k}AA^{T}\mathbf{r}_{k},
	\end{equation}
	and hence
	\begin{equation}\label{equ:dk}
	\mathbf{d}_{k-1}=\mathbf{r}_{k}-\mathbf{r}_{k-1}=-\alpha_{k-1}AA^{T}\mathbf{r}_{k-1}.
	\end{equation}
	Let $B=AA^{T}$. Dropping the common positive factor $\alpha_{k-1}^{2}$, from \eqref{equ:dk}, we obtain
	\begin{equation*}
	\alpha_{k}^{\mathrm{RD}}
	=\frac{\mathbf{d}_{k-1}^{T}\mathbf{d}_{k-1}}{\|A^{T}\mathbf{d}_{k-1}\|_2^{2}}
	=\frac{\mathbf{r}_{k-1}^{T}B^{2}\mathbf{r}_{k-1}}{\mathbf{r}_{k-1}^{T}B^{3}\mathbf{r}_{k-1}},
	\qquad
	\alpha_{k}^{\mathrm{DR}}
	=\frac{\mathbf{r}_{k-1}^{T}\mathbf{r}_{k-1}}{\|A^{T}\mathbf{r}_{k-1}\|_2^{2}}
	=\frac{\mathbf{r}_{k-1}^{T}\mathbf{r}_{k-1}}{\mathbf{r}_{k-1}^{T}B\,\mathbf{r}_{k-1}}.
	\end{equation*}
	By applying the Cauchy–Schwarz inequality twice, we obtain $\alpha_{k}^{\mathrm{RD}}\le\alpha_{k}^{\mathrm{DR}}$, which completes the proof.
\end{proof}

From the definitions of the two step sizes and Theorem 3.1, it follows that the new step sizes $\alpha_k^{\mathrm{DR}}$ and $\alpha_k^{\mathrm{RD}}$ correspond, respectively, to the BB step sizes $\alpha_k^{\mathrm{BB1}}$ and $\alpha_k^{\mathrm{BB2}}$ in convex quadratic optimization. In what follows, we present a convergence analysis of the gradient descent method equipped with the step size $\alpha_k^{\mathrm{DR}}$, for consistent linear systems.

\subsection{Convergence analysis}
Let
\begin{equation*}
A_0:=\begin{bmatrix}\sigma_1&0\\0&\sigma_2\\0&0\end{bmatrix}
\end{equation*}
denote the canonical form of $A$, with $\sigma_1\ge\sigma_2>0$. Since $A=\widetilde{U}A_0V^{T}$ by the singular value decomposition and the gradient descent iteration is equivariant under orthogonal transformations, we may assume without loss of generality that $A=A_0$, with $\mathbf{b}$ replaced by $\widetilde{U}^{T}\mathbf{b}$. In the consistent case, the third component of $\widetilde{U}^{T}\mathbf{b}$ vanishes, and the zero row plays no role in the iteration.

Consequently, the third component of the residual vanishes identically, and it suffices to track the first two components of $\mathbf{r}_k$. Set
\begin{equation}\label{equ:epsilon}
\varepsilon_{k-1}:=\frac{\big(\mathbf{r}_{k-1}^{(1)}\big)^{2}}{\big(\mathbf{r}_{k-1}^{(2)}\big)^{2}}>0.	
\end{equation}
By the definition of $\alpha_k^{\mathrm{DR}}$,
\begin{equation}
\alpha_k^{\mathrm{DR}}=\frac{\|\mathbf{r}_{k-1}\|_2^{2}}{\|A_0^{T}\mathbf{r}_{k-1}\|_2^{2}}=\frac{\varepsilon_{k-1}+1}{\sigma_1^{2}\varepsilon_{k-1}+\sigma_2^{2}}.
\end{equation}
From \eqref{equ:rrk} and \eqref{equ:iterate},
\begin{equation}
\mathbf{r}_{k+1}=\big(I-\alpha_k^{\mathrm{DR}}AA^{T}\big)\mathbf{r}_k.
\end{equation}
In component form,
\begin{equation}\label{equ:rh}
\begin{cases}
\mathbf{r}_{k+1}^{(1)}=\eta_k\mathbf{r}_k^{(1)},\\
\mathbf{r}_{k+1}^{(2)}=\xi_k\mathbf{r}_k^{(2)},
\end{cases}
\end{equation}
where
\begin{equation}\label{equ:ex}
\eta_k=\frac{\sigma_2^{2}-\sigma_1^{2}}{\sigma_1^{2}\varepsilon_{k-1}+\sigma_2^{2}},\qquad
\xi_k=\frac{(\sigma_1^{2}-\sigma_2^{2})\varepsilon_{k-1}}{\sigma_1^{2}\varepsilon_{k-1}+\sigma_2^{2}}\in(0,1).
\end{equation}
Since \(\|\mathbf{r}_{k+1}\|_2^{2}=(\mathbf{r}_{k+1}^{(2)})^{2}(1+\varepsilon_{k+1})\) and \(\|\mathbf{r}_k\|_2^{2}=(\mathbf{r}_k^{(2)})^{2}(1+\varepsilon_k)\), we have
\begin{equation}
\frac{\|\mathbf{r}_{k+1}\|_2^{2}}{\|\mathbf{r}_k\|_2^{2}}=\xi_k^{2}\frac{1+\varepsilon_{k+1}}{1+\varepsilon_k}.
\end{equation}
Iterating this relation yields
\begin{equation}\label{equ:recurise}
\frac{\|\mathbf{r}_{k+N}\|_2^{2}}{\|\mathbf{r}_k\|_2^{2}}=\Big(\prod_{j=0}^{N-1}\xi_{k+j}\Big)^{2}\frac{1+\varepsilon_{k+N}}{1+\varepsilon_k}.
\end{equation}
Hence it remains to prove that there exists some \(N\ge1\) such that \(\varepsilon_{k+N}\le\varepsilon_k\).

It follows from the definition \eqref{equ:epsilon} of $\varepsilon_{k}$, \eqref{equ:rh} and \eqref{equ:ex} that 
\begin{equation}\label{equ:dy}
	\varepsilon_{k+1}=\frac{\varepsilon_{k}}{\varepsilon_{k-1}^{2}}.
\end{equation}
Taking logarithms in \eqref{equ:dy}, i.e., setting \(\varphi_k=\ln\varepsilon_k\), yields a second-order homogeneous linear difference equation
\(\varphi_{k+1}-\varphi_k+2\varphi_{k-1}=0,\)
whose characteristic roots are \(\dfrac{1\pm i\sqrt7}{2}=\sqrt2\,e^{\pm i\theta}\) with \(\theta=\arccos\dfrac{1}{\sqrt8}\). Hence the general solution can be written uniformly as
\begin{equation}\label{equ:phik}
\varphi_k=(\sqrt2)^{k}\,C\cos(k\theta-\delta),
\end{equation}
where \(C=\sqrt{C_1^{2}+C_2^{2}}\) is the amplitude and \(\delta=\arctan(C_2/C_1)\) the phase, both being determined by the initial values \(\varphi_0\) and \(\varphi_1\). Consequently,
\begin{equation}\label{equ:ep}
\varepsilon_k=e^{(\sqrt2)^{k}C\cos(k\theta-\delta)}.
\end{equation}

Next, based on the explicit expression of \(\varepsilon_k\), we prove the convergence of the DR method.

\begin{lemma}\label{lem:irrational}
	Let $\theta=\arccos\bigl(1/(2\sqrt2)\bigr)$. Then $\theta/\pi$ is irrational.
\end{lemma}

\begin{proof}
	We have $2\cos\theta=1/\sqrt2$. If $\theta/\pi$ were rational, then by Niven's theorem
	 $2\cos\theta$ would be an algebraic integer \cite{Niven1956}. However, the minimal polynomial of
	$1/\sqrt2$ is $2x^{2}-1=0$, which is not monic; hence $1/\sqrt2$ is not an algebraic integer,
	a contradiction.
\end{proof}

\begin{lemma}\label{lem:density}
	For every $\beta\in\mathbb{R}$ and every $t\in[-1,1]$, there exist infinitely many integers
	$N\ge1$ such that $\cos(\beta+N\theta)\le t$.
\end{lemma}

\begin{proof}
	By Lemma~\ref{lem:irrational}, $\theta/\pi\notin\mathbb{Q}$, so the rotation
	$\{N\theta\bmod 2\pi\}_{N\ge1}$ is dense on the circle; the claim follows immediately.
\end{proof}

\begin{lemma}\label{lem:eps}
	Let $\varepsilon_k=e^{(\sqrt2)^{k}C\cos(k\theta-\delta)}$ with $C\ge0$. Then for every $k$ there
	exist infinitely many integers $N\ge1$ such that $\varepsilon_{k+N}\le\varepsilon_k$.
\end{lemma}

\begin{proof}
	If $C=0$, then $\varepsilon_k\equiv1$ and the claim is trivial. Assume $C>0$ and set
	$\beta:=k\theta-\delta$ and $\varphi_j:=(\sqrt2)^{j}C\cos(j\theta-\delta)$.
	If $\cos\beta\ge0$, Lemma~\ref{lem:density} with $t=0$ yields infinitely many $N$ with
	$\cos(\beta+N\theta)\le0$, so $\varphi_{k+N}\le0\le\varphi_k$.
	If $\cos\beta<0$, Lemma~\ref{lem:density} with $t=-1/2$ yields infinitely many $N$ with
	$\cos(\beta+N\theta)<-1/2$; along this subsequence,
	$\varphi_{k+N}<-\tfrac12(\sqrt2)^{k+N}C\to-\infty$ as $N\to\infty$, so
	$\varphi_{k+N}\le\varphi_k$ for all sufficiently large such $N$.
	In both cases, $\varepsilon_{k+N}=e^{\varphi_{k+N}}\le e^{\varphi_k}=\varepsilon_k$.
\end{proof}

\begin{theorem}[Residual contraction over a window of $N$ steps]\label{prop:contract}
	Let $\sigma_1\ge\sigma_2>0$ and consider the canonical system $A_0\mathbf{x}=\mathbf{b}$, where
	$A_0=[\operatorname{diag}(\sigma_1,\sigma_2);\,0]$ and $b_3=0$. Let $\{\mathbf{x}_k\}$ be generated
	by gradient descent $\mathbf{x}_{j+1}=\mathbf{x}_j-\alpha_j^{\mathrm{DR}}A_0^{T}\mathbf{r}_j$ with
	residuals $\mathbf{r}_j=\mathbf{b}-A_0\mathbf{x}_j$, and assume
	$\mathbf{r}_{k-1}^{(2)}\neq 0$, so that the quantity $\varepsilon_{k-1}$ in \eqref{equ:epsilon} is well defined.
	Then there exists an integer $N>1$, depending on $k$, such that
	$\|\mathbf{r}_{k+N}\|_2<\Big(\prod_{j=0}^{N-1}\xi_{k+j}\Big)^{2}\|\mathbf{r}_k\|_2$.
\end{theorem}

\begin{proof}
It follows from \eqref{equ:recurise} and Lemma \ref{lem:eps} that there exists some \(N\ge1\) such that \(\|\mathbf{r}_{k+N}\|_2<\Big(\prod_{j=0}^{N-1}\xi_{k+j}\Big)^{2}\|\mathbf{r}_k\|_2\). This completes the proof.
\end{proof}

\section{Numerical experiments}

In this section we report numerical experiments\footnote{All numerical experiments in this section are performed in MATLAB (version R2024a).} that verify the effectiveness of the
two proposed step sizes, the delayed residual step size $\alpha_k^{\mathrm{DR}}$ and
the residual-difference step size $\alpha_k^{\mathrm{RD}}$ defined in \eqref{equ:RDDR}.
The proposed methods are compared with the classical cyclic Kaczmarz method
(KACZ), the randomized Kaczmarz method (RK), and the residual-based surrogate
hyperplane Kaczmarz method (RSHK), which is equivalent to the method of steepest
descent \cite{Wang2023}. Two groups of experiments are performed: in
Experiment~1 we record the convergence curves of the relative residual norm on
ill-conditioned consistent linear systems with condition numbers
$\kappa=10^{2},10^{3},10^{4}$; in Experiment~2 we reproduce the Kaczmarz anomaly
phenomenon of Dax \cite{Achiya2023} and demonstrate that DR/RD substantially
alleviate the slow convergence of KACZ as the number of rows approaches the
number of columns.

\subsection{Experimental setup}\label{sec:setup}

Throughout this section the following conventions are adopted.

\begin{itemize}
  \item \textbf{Stopping criterion.} All methods start from the zero vector
  $\mathbf{x}_{0}=\mathbf{0}$ and are stopped once the relative residual satisfies
  \begin{equation}
    \frac{\|\mathbf{r}_{k}\|_{2}}{\|\mathbf{b}\|_{2}}<10^{-6},
    \qquad \mathbf{r}_{k}=\mathbf{b}-A\mathbf{x}_{k},
  \end{equation}
  or once the number of iterations reaches the cap
  $k_{\max}=2\times10^{6}$; in the latter case the run is recorded as
  \emph{not converged} and the final relative residual is reported.
  \item \textbf{Initialization of DR/RD.} Since $\alpha_k^{\mathrm{DR}}$ and
  $\alpha_k^{\mathrm{RD}}$ in \eqref{equ:RDDR} require the information of the
  previous step, at $k=0$ the RSHK step size
  $\alpha_0^{\mathrm{RSHK}}=\|\mathbf{r}_{0}\|_{2}^{2}/\|A^{T}\mathbf{r}_{0}\|_{2}^{2}$
  is used to start the recursion.
  \item \textbf{Reproducibility.} All random quantities (the test matrices, the
  true solutions, and the row sampling of RK) are generated with a fixed random
  seed, so that all reported results are reproducible.
  \item \textbf{Numerical safeguard.} A denominator smaller than $10^{-30}$ in
  the step-size formulas is treated as numerical convergence and the iteration
  is terminated, in order to avoid overflow.
\end{itemize}

All test systems are consistent, i.e.\ $\mathbf{b}=A\mathbf{x}_{\ast}$ with a
randomly generated $\mathbf{x}_{\ast}\in\mathbb{R}^{n}$, so that the relative
residual converges to zero for every method.

\subsection{Experiment 1: ill-conditioned consistent systems}\label{sec:exp1}

\subsubsection{Problem description.}
We take square systems with $m=n=100$. The coefficient matrix is constructed as
$A=U\,\mathrm{diag}(\sigma_{1},\ldots,\sigma_{n})V^{T}$, where
$U\in\mathbb{R}^{m\times n}$ and $V\in\mathbb{R}^{n\times n}$ are random matrices
with orthonormal columns, and the singular values are geometrically spaced
between $\sigma_{1}=1$ and $\sigma_{n}=1/\kappa$, so that
$\mathrm{cond}(A)=\kappa$ holds exactly. The three values
$\kappa=10^{2},10^{3},10^{4}$ are tested.

\subsubsection{Results.}
\begin{figure}[htbp]
\centering
\begin{subfigure}{0.32\linewidth}
  \includegraphics[width=\linewidth]{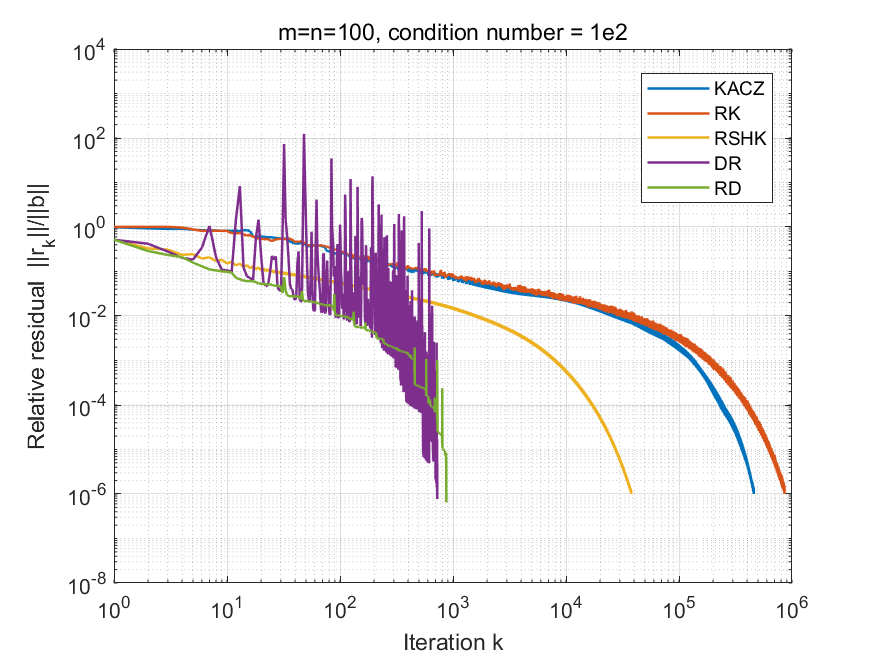}
  \caption{$\kappa=10^{2}$}
\end{subfigure}\hfill
\begin{subfigure}{0.32\linewidth}
  \includegraphics[width=\linewidth]{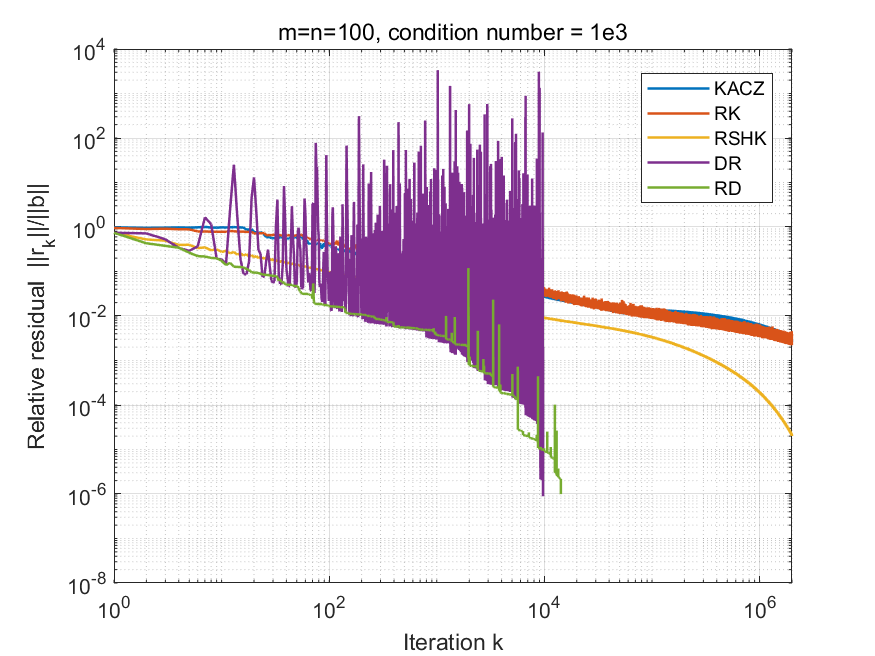}
  \caption{$\kappa=10^{3}$}
\end{subfigure}\hfill
\begin{subfigure}{0.32\linewidth}
  \includegraphics[width=\linewidth]{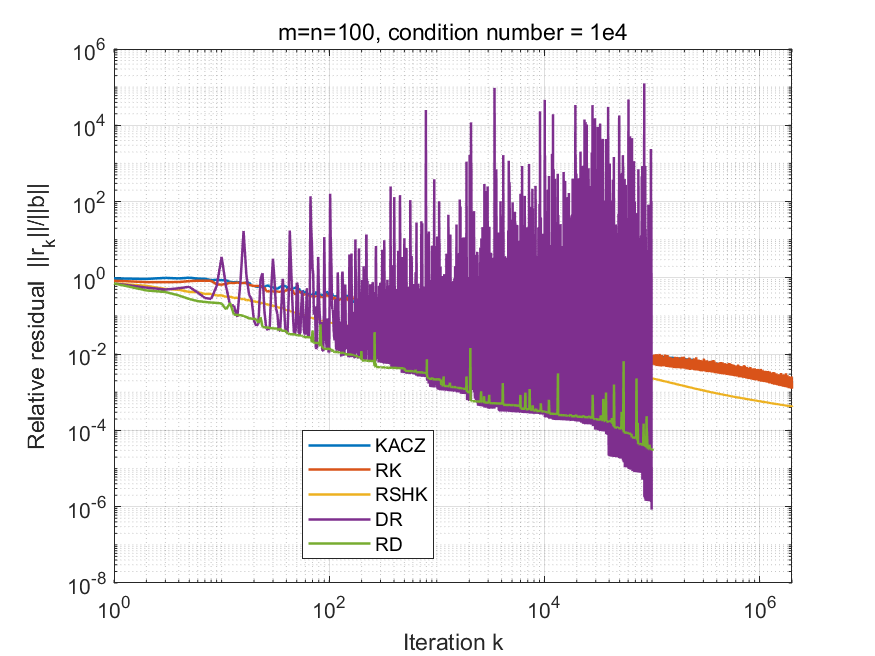}
  \caption{$\kappa=10^{4}$}
\end{subfigure}
\caption{Experiment 1 ($m=n=100$): relative residual norms
$\|\mathbf{r}_{k}\|_{2}/\|\mathbf{b}\|_{2}$ versus the number of iterations $k$
(log--log scale) for KACZ, RK, RSHK, DR and RD.}
\label{fig:exp1}
\end{figure}

Table~\ref{tab:exp1} reports the numbers of iterations required by the five
methods to reach the tolerance $10^{-6}$; for the non-converged cases the final
relative residual at the cap $k_{\max}=2\times10^{6}$ is given in parentheses.
Figure~\ref{fig:exp1} displays the corresponding convergence curves of the
relative residual norm in log--log scale.

\begin{table}[htbp]
\centering
\small
\caption{Experiment 1 ($m=n=100$): numbers of iterations to reach
$\|\mathbf{r}_{k}\|_{2}/\|\mathbf{b}\|_{2}<10^{-6}$. For the non-converged
cases (marked $>2\times10^{6}$), the value in parentheses is the final relative
residual at the iteration cap.}
\label{tab:exp1}
\begin{tabular}{lccccc}
\toprule
$\kappa$ & KACZ & RK & RSHK & DR & RD \\
\midrule
$10^{2}$ & $458{,}794$ & $858{,}498$ & $37{,}894$ & $\mathbf{726}$ & $\mathbf{873}$ \\
$10^{3}$ & $>2\times10^{6}\ (2.6\times10^{-3})$ & $>2\times10^{6}\ (3.5\times10^{-3})$ &
$>2\times10^{6}\ (2.0\times10^{-5})$ & $\mathbf{9{,}712}$ & $\mathbf{14{,}219}$ \\
$10^{4}$ & $>2\times10^{6}\ (1.8\times10^{-3})$ & $>2\times10^{6}\ (1.7\times10^{-3})$ &
$>2\times10^{6}\ (4.4\times10^{-4})$ & $\mathbf{99{,}114}$ &
$>2\times10^{6}\ (3.1\times10^{-5})$ \\
\bottomrule
\end{tabular}
\end{table}

\subsubsection{Analysis.}
The results confirm the conclusions of Sections~2 and~3.

\begin{enumerate}
  \item \textbf{DR/RD are markedly more robust to ill-conditioning.} For
  $\kappa=10^{2}$, DR and RD converge in $726$ and $873$ iterations,
  respectively, which is about $52$ times faster than RSHK (the method of
  steepest descent, $37{,}894$ iterations) and about $630$ times faster than
  KACZ ($458{,}794$ iterations). For $\kappa=10^{3}$, KACZ, RK and RSHK all fail
  to converge within the cap $2\times10^{6}$ (their final relative residuals
  only reach $10^{-3}\sim10^{-5}$), whereas DR and RD still converge in
  $9{,}712$ and $14{,}219$ iterations. For $\kappa=10^{4}$, DR remains
  convergent ($99{,}114$ iterations) while the three comparison methods are
  essentially stalled. This is consistent with the known acceleration effect of
  BB-type step sizes on ill-conditioned quadratics \cite{barzilai1988,
  dai2002} and with the analysis of Section~3.
  \item \textbf{Non-monotonicity is visible in the curves.} In
  Figure~\ref{fig:exp1} the residual norms of DR and RD first increase, then
  oscillate and eventually drop sharply, which is the non-monotone
  \emph{window-contraction} behavior established in
  Theorem~\ref{prop:contract}: despite transient increases, the residual norm
  decreases after finitely many steps.
  \item \textbf{One phenomenon deserves an honest comment.} At
  $\kappa=10^{4}$, RD does not reach $10^{-6}$ within the cap (final relative
  residual $3.1\times10^{-5}$). A plausible explanation is that
  $\alpha_k^{\mathrm{RD}}\le\alpha_k^{\mathrm{DR}}$ (Theorem~3.1): on a
  geometrically graded singular spectrum the BB2-type step size tends to be
  smaller, and the residual-difference construction is more sensitive to
  round-off in the late stage of convergence. This suggests that DR and RD play
  complementary roles: DR is more robust for extremely ill-conditioned
  problems, while RD is faster in the regime $m\approx n$ considered next.
\end{enumerate}

\subsection{Experiment 2: the Kaczmarz anomaly and its alleviation}\label{sec:exp2}

\subsubsection{Problem description.}
Following Dax \cite{Achiya2023}, we fix the number of columns $n=100$ and let the
number of rows $m$ increase from $5$ to $200$:
\[
  m\in\{5,10,20,40,60,80,90,95,98,100,105,110,120,150,200\}.
\]
The matrices $A=\mathrm{randn}(m,n)$ are random Gaussian matrices, so their
condition numbers are natural (not artificially controlled), which reproduces
the regime in which the anomaly occurs. The system is again consistent,
$\mathbf{b}=A\mathbf{x}_{\ast}$. Note that for $m<n$ the system is
underdetermined; starting from $\mathbf{x}_{0}=\mathbf{0}$, all methods
converge to the least-norm solution $A^{\dagger}\mathbf{b}$, and the relative
residual criterion remains valid.

\subsubsection{Results.}
\begin{figure}[htbp]
\centering
\includegraphics[width=0.72\linewidth]{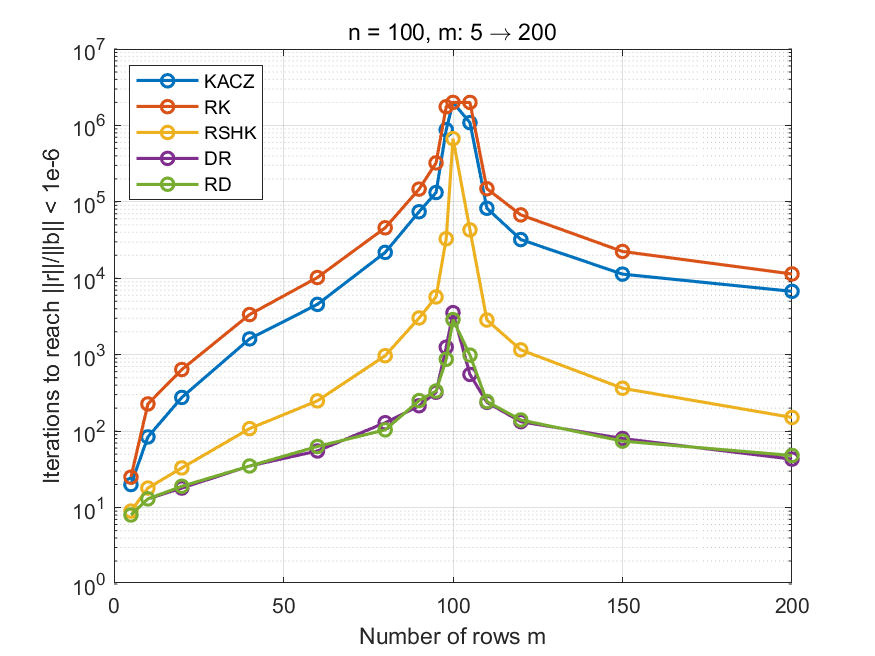}
\caption{Experiment 2 ($n=100$): numbers of iterations required by KACZ, RK,
RSHK, DR and RD to reach $\|\mathbf{r}_{k}\|_{2}/\|\mathbf{b}\|_{2}<10^{-6}$
versus the number of rows $m$ (logarithmic vertical axis; non-converged cases
are plotted at the cap $2\times10^{6}$).}
\label{fig:exp2}
\end{figure}

Table~\ref{tab:exp2} reports the numbers of iterations required by each method
for every tested $m$, and Figure~\ref{fig:exp2} plots them on a logarithmic
vertical axis (non-converged cases are plotted at the cap $2\times10^{6}$).

\begin{table}[htbp]
\centering
\small
\caption{Experiment 2 ($n=100$, $A=\mathrm{randn}(m,100)$): numbers of
iterations to reach $\|\mathbf{r}_{k}\|_{2}/\|\mathbf{b}\|_{2}<10^{-6}$
(the entries $>2\times10^{6}$ mean that the method did not converge within the
iteration cap).}
\label{tab:exp2}
\begin{tabular}{lccccc}
\toprule
$m$ & KACZ & RK & RSHK & DR & RD \\
\midrule
5    & $20$          & $25$          & $9$           & $8$           & $8$ \\
10   & $84$          & $227$         & $18$          & $13$          & $13$ \\
20   & $276$         & $640$         & $33$          & $18$          & $19$ \\
40   & $1{,}612$     & $3{,}355$     & $108$         & $35$          & $35$ \\
60   & $4{,}557$     & $10{,}236$    & $249$         & $55$          & $63$ \\
80   & $21{,}822$    & $45{,}776$    & $971$         & $129$         & $104$ \\
90   & $74{,}264$    & $147{,}016$   & $3{,}027$     & $215$         & $251$ \\
95   & $132{,}484$   & $323{,}402$   & $5{,}703$     & $323$         & $335$ \\
98   & $880{,}612$   & $1{,}759{,}101$ & $32{,}919$  & $1{,}256$     & $871$ \\
100  & $>2\times10^{6}$ & $>2\times10^{6}$ & $671{,}843$ & $\mathbf{3{,}552}$ & $\mathbf{2{,}879}$ \\
105  & $1{,}086{,}826$ & $>2\times10^{6}$ & $42{,}962$ & $\mathbf{552}$ & $\mathbf{988}$ \\
110  & $82{,}458$    & $148{,}668$   & $2{,}843$     & $238$         & $244$ \\
120  & $32{,}107$    & $67{,}542$    & $1{,}157$     & $133$         & $140$ \\
150  & $11{,}333$    & $22{,}430$    & $364$         & $80$          & $74$ \\
200  & $6{,}747$     & $11{,}386$    & $151$         & $43$          & $48$ \\
\bottomrule
\end{tabular}
\end{table}

\subsubsection{Analysis.}
\begin{enumerate}
  \item \textbf{The Kaczmarz anomaly is fully reproduced.} The iteration count
  of KACZ first increases steadily (from $20$ at $m=5$), then surges
  dramatically as $m$ approaches $n$ ($880{,}612$ at $m=98$), exceeds the cap
  at $m=100$, and finally recovers once $m$ passes $n$ ($6{,}747$ at $m=200$).
  This ``deterioration--recovery'' behavior is in full qualitative agreement
  with Table~1 of Dax \cite{Achiya2023} (whose counts are smaller because a
  different stopping criterion and different random matrices are used there).
  \item \textbf{DR/RD almost completely eliminate the anomaly.} Throughout the
  whole interval $5\le m\le200$, the iteration counts of DR and RD stay in
  $8\sim3{,}552$ and $8\sim2{,}879$, respectively, i.e.\ they remain almost
  flat and are essentially unaffected by the anomalous peak. At $m=100$, DR and
  RD are more than $563$ and $695$ times faster than KACZ, and about $189$ and
  $233$ times faster than RSHK ($671{,}843$ iterations); at $m=105$, KACZ still
  needs $1{,}086{,}826$ iterations while DR needs only $552$ (about $1{,}970$
  times faster). In the overdetermined regime ($m=200$) DR and RD converge in
  $43$ and $48$ iterations, still clearly better than KACZ ($6{,}747$) and
  RSHK ($151$).
  \item \textbf{RSHK itself still suffers from the anomaly.} Although RSHK is
  much faster than KACZ, its iteration count still grows sharply as $m\to n$
  ($671{,}843$ at $m=100$). This shows that the optimal step size of steepest
  descent alone cannot fully escape the ill-conditioned regime, and that the
  second-order dynamics induced by the delayed/difference step sizes of the BB
  type is essential---which is exactly the mechanism exploited by DR and RD in
  Section~3.
  \item \textbf{Complementary behavior of DR and RD.} For $m\ge80$ RD is
  slightly faster than DR (e.g.\ $871$ vs.\ $1{,}256$ at $m=98$), in contrast
  with Experiment~1 at $\kappa=10^{4}$. Thus RD is more efficient under the
  natural ill-conditioning caused by $m\approx n$, while DR is more robust on
  extremely ill-conditioned systems with a graded singular spectrum.
\end{enumerate}

In summary, the numerical results show that the proposed DR and RD step sizes
reduce the number of iterations needed for the convergence precision $10^{-6}$
by two to three orders of magnitude on both the ill-conditioned systems of
Experiment~1 and the full anomaly interval $5\le m\le200$ of Experiment~2,
thereby substantially alleviating the slow convergence of Kaczmarz-type
methods.

\section{Conclusions and outlook}\label{sec:conclusion}

In this paper, we have proposed two-point residual step size gradient methods
for solving consistent linear systems, based respectively on a
\emph{delayed residual} (DR) step size and a \emph{residual-difference} (RD)
step size, which correspond, respectively, to the two Barzilai-Borwein (BB)
step sizes in convex quadratic optimization. Like the BB method, the proposed
step sizes are non-monotone: although the relative residual norm may increase
in individual steps, we proved, in the canonical two-singular-value model of
Section~3.2, that it decreases after finitely many steps
(Theorem~\ref{prop:contract}). Numerical experiments on ill-conditioned linear
systems with $\kappa=10^2,10^3,10^4$ demonstrated that both new methods
converge far more rapidly than the cyclic Kaczmarz method, the randomized
Kaczmarz method, and the method of steepest descent, and that they
substantially alleviate the slow-convergence regime of Kaczmarz's anomaly over
the whole interval $5\le m\le200$.

Following the extensive improvements of the BB method
\cite{Zhou2006GradientMethodsAdaptive,Frassoldati2008Newadaptivestepsize,Friedlander1998GradientMethodRetards,An2025RegularizedBarzilaiBorwein},
several refined strategies can be transplanted in a natural way to the proposed
DR and RD step sizes, including the alternating BB (ABB) method, the ABBmin method, and regularized BB (RBB) methods.
Investigating the application of these strategies to the efficient solution of
linear systems---particularly in the ill-conditioned and anomalous regimes
identified in this paper---is a promising direction for future research.



\begin{thebibliography}{99}
	
\bibitem{Achiya2023}
A.~Dax, \emph{Kaczmarz’s anomaly: A surprising feature of 
	Kaczmarz’s method}, Linear Algebra Appl.\ \textbf{662} (2023), 136--162.

\bibitem{An2025RegularizedBarzilaiBorwein} 
C.~P.~An, X.~Xu, \emph{Regularized Barzilai-Borwein method}, Numerical Algorithms,\  \textbf{102} (2025), 25--52. 

\bibitem{Friedlander1998GradientMethodRetards}
A. Friedlander, J.M. Mart\'{i}nez, B. Molina, M. Raydan, \emph{Gradient Method with Retards and Generalizations}, SIAM J. Numer. Anal., \textbf{36}(1)(1998), pp. 275--289.

\bibitem{Frassoldati2008Newadaptivestepsize}
G. Frassoldati, L. Zanni, G. Zanghirati,
\emph{New adaptive stepsize selections in gradient methods},
J. Ind. Manag. Optim. 4(2008), pp. 299-312.
	
\bibitem{kacz1937}	
S.~Kaczmarz, \emph{Angen\"{a}herte aufl\"{o}sung von systemen linearer gleichungen}, Bull. Int.
Acad. Polon. Sci. Lett. Ser. A, (1937), 335--357.

\bibitem{strohmer2009}
T.~Strohmer and R.~Vershynin, \emph{A randomized Kaczmarz algorithm with exponential convergence},
J.~Fourier Anal.\ Appl.\ \textbf{15} (2009), 262--278.

\bibitem{needell2014}
D.~Needell and J.~A.~Tropp, \emph{Paved with good intentions: Analysis of a randomized block Kaczmarz
method}, Linear Algebra Appl.\ \textbf{441} (2014), 199--221.

\bibitem{Niven1956}
I.~Niven, \emph{Irrational Numbers}, Mathematical Association of America, Carus Mathematical Monographs, New York, 1956.

\bibitem{han2021}
R.~Han, L.~Li, P.~Yang, F.~Zhang, and X.~Gao, \emph{A novel constrained reconstruction model towards high-resolution subtomogram averaging}, Bioinformatics\ \textbf{37} (2021), 1616--1623.

\bibitem{Huang2021EquippingBarzilaiBorwein} Y.K. Huang, Y.H. Dai, X.W. Liu, Equipping the Barzilai-Borwein Method with the Two Dimensional Quadratic Termination Property, SIAM J. Optim. \textbf{3}(4) (2021) 3068--3096

\bibitem{Raydan1997BarzilaiBorweinGradienta} M. Raydan, The Barzilai and Borwein Gradient Method for the Large Scale Unconstrained Minimization Problem, SIAM J. Optim. \textbf{7}(1) (1997) 26--33

\bibitem{saad2011}
Y.~Saad, \emph{Iterative Methods for Sparse Linear Systems}, 2nd ed., SIAM, Philadelphia, 2003.

\bibitem{Wang2023}
Z.~Wang and J.~F.~Yin, \emph{A surrogate hyperplane Kaczmarz method for solving consistent linear equations}, Applied Mathematics Letters\ \textbf{144} (2023), 108704.

\bibitem{gower2015}
R.~M.~Gower and P.~Richt\'arik, \emph{Randomized iterative methods for linear systems}, SIAM J.~Matrix Anal.\ Appl.\ \textbf{36} (2015), 1660--1690.

\bibitem{du2020}
K.~Du, W.~Si, and X.~Sun, \emph{Pseudoinverse-free randomized extended block Kaczmarz for solving least squares}, SIAM J.~Sci.\ Comput.\ \textbf{42} (2020), A130--A150.

\bibitem{bai2018}
Z.-Z.~Bai and W.-T.~Wu, \emph{On greedy randomized Kaczmarz method for solving large sparse linear systems}, SIAM J.~Sci.\ Comput.\ \textbf{40} (2018), A592--A606.

\bibitem{liu2016}
J.~Liu and S.~J.~Wright, \emph{An accelerated randomized Kaczmarz algorithm}, Math.\ Comp.\ \textbf{85} (2016), 153--178.

\bibitem{chi2016}
Y.~Chi and Y.~M.~Lu, \emph{Kaczmarz method for solving quadratic equations}, IEEE Trans.\ Inf.\ Theory \textbf{62} (2016), 5477--5489.

\bibitem{gordon1970}
R.~Gordon, R.~Bender, and G.~T.~Herman, \emph{Algebraic reconstruction techniques (ART) for three-dimensional electron microscopy and X-ray photography}, J.~Theoret.\ Biol.\ \textbf{29} (1970), 471--481.

\bibitem{barzilai1988}
J.~Barzilai and J.~M.~Borwein, \emph{Two-point step size gradient methods}, IMA J.~Numer.\ Anal.\ \textbf{8} (1988), 141--148.

\bibitem{dai2002}
Y.-H.~Dai and L.-Z.~Liao, \emph{$R$-linear convergence of the Barzilai and Borwein gradient method}, IMA J.~Numer.\ Anal.\ \textbf{22} (2002), 1--10.

\bibitem{Xu2025PBB} 
X. Xu, A parameterized Barzilai-Borwein method via
interpolated least squares, Journal of Industrial and Management Optimization. \textbf{21}(11) (2025)  6246-6269

\bibitem{Xu2025OnRBB}
X. Xu, \emph{On Convergence of Regularized Barzilai-Borwein Method}, arXiv preprint (2025), arXiv:2512.21972.

\bibitem{Zhou2006GradientMethodsAdaptive} 
B. Zhou, L. Gao, Y.H. Dai, Gradient Methods with Adaptive Step-Sizes, Comput Optim Appl. \textbf{35}(1) (2006) 69-86

\end{thebibliography}
\end{document}